\documentclass{ifacconf}

\usepackage{microtype}
\usepackage{graphicx}      
\usepackage{natbib}        
\usepackage{amsmath}       
\usepackage{amssymb}       

\usepackage{mathtools}
\usepackage{graphicx}
\usepackage{physics}

\usepackage{algorithm}
\usepackage{algpseudocode}

\usepackage{spverbatim}

\usepackage{mathtools} 

\usepackage{xcolor}

\usepackage{pgf}

\newcommand*{\SectorRadius}{0.9ex}
\newcommand*{\SectorHalfAngle}{27}
\newcommand*{\SectorLineWidth}{0.5pt}

\newcommand*{\sector}{%
  \hspace{0.5pt}\begin{pgfpicture}
    \pgfpathmoveto{\pgforigin}%
    \pgfpathlineto{\pgfpointpolar{90-(\SectorHalfAngle)}{\SectorRadius}}%
    \pgfarc{90-(\SectorHalfAngle)}{90+\SectorHalfAngle}{\SectorRadius}%
    \pgfpathclose
    \pgfsetlinewidth{\SectorLineWidth}%
    \pgfusepath{stroke}%
  \end{pgfpicture}%
}

\usepackage{accents}

\makeatletter
\let\old@ssect\@ssect 
\makeatother

\usepackage[
colorlinks=true
]{hyperref}

\definecolor{dgreen}{rgb}{0, 0.5, 0}

\hypersetup{
citecolor={dgreen}
}

\makeatletter
\def\@ssect#1#2#3#4#5#6{%
  \NR@gettitle{#6}
  \old@ssect{#1}{#2}{#3}{#4}{#5}{#6}
}
\makeatother
\newcommand{\behcetacikmese}{Beh\c{c}et~A\c{c}\i kme\c{s}e}

\newcommand{\pipg}{\textsc{pipg}}
\newcommand{\seco}{\textsc{s}\begin{footnotesize}e\end{footnotesize}\textsc{co}}
\newcommand{\ecos}{\textsc{ecos}}
\newcommand{\ptr}{\textsc{ptr}}
\newcommand{\gurobi}{\textsc{gurobi}}
\newcommand{\mosek}{\textsc{mosek}}

\newcommand{\R}{\mathbb{R}}

\newcommand{\D}{\mathbb{D}}

\newcommand{\defeq}{\vcentcolon=}

\newcommand{\f}[2]{#1\!\left(#2\right)} 

\definecolor{darts}{HTML}{009670}
\definecolor{bricks}{HTML}{E74C3C}
\definecolor{fista}{HTML}{075187}
\definecolor{fistagray}{HTML}{929598}
\definecolor{steelblue}{HTML}{4682b4}
\definecolor{goldenrod}{HTML}{daa520}
\definecolor{purple}{HTML}{800080}
\definecolor{light}{HTML}{d3d3d3}

\newcommand{\range}[2]{#1\!:\!#2} 

\newcommand{\circdot}[1]{\accentset{\circ}{#1}} 

\newcommand{\I}{\mathcal{I}}
\newcommand{\B}{\mathcal{B}}

\begin{document}

\begin{frontmatter}
\title{Real-Time Sequential Conic Optimization for Multi-Phase Rocket Landing Guidance\!\thanksref{footnoteinfo}}

\thanks[footnoteinfo]{\!\!This work was supported by NASA grant NNX17AH02A.}

\author[UW]{Abhinav G.\@ Kamath\!}
\author[UW]{~Purnanand Elango\!}
\author[UT]{~Yue Yu\!}
\author[UW]{Skye Mceowen\!}
\author[UW]{~Govind M.\@ Chari\!}
\author[NASA]{~John M.\@ Carson III\!}
\author[UW]{\behcetacikmese\!}

\address[UW]{William E.\@ Boeing Department of Aeronautics \& Astronautics, University of Washington, Seattle, WA 98195, USA \\\normalfont{(e-mail: \texttt{\{agkamath,\,pelango,\,skye95,\,gchari,\,behcet\}@uw.edu})}.}
\address[UT]{Oden Institute for Computational Engineering and Sciences, The University of Texas at Austin, Austin, TX 78712, USA \\\normalfont{(e-mail: \normalfont{\texttt{yueyu@utexas.edu}})}.}
\address[NASA]{NASA Johnson Space Center, Houston, TX 77058, USA \\\normalfont{(e-mail: \normalfont{\texttt{john.m.carson@nasa.gov}})}.}
\begin{abstract}
We introduce a multi-phase rocket landing guidance framework that can handle nonlinear dynamics and does not mandate any additional mixed-integer or nonconvex constraints to handle discrete temporal events/switching. To achieve this, we first introduce \textit{sequential conic optimization} (\seco), a new paradigm for solving nonconvex optimal control problems that is entirely devoid of matrix factorizations and inversions. This framework combines sequential convex programming (SCP) and first-order conic optimization and can solve unified multi-phase trajectory optimization problems in real-time. The novel features of this framework are: (1) \textit{time-interval dilation}, which enables multi-phase trajectory optimization with free-transition-time; (2) \textit{single-crossing compound state-triggered constraints}, which are entirely convex if the trigger and constraint conditions are convex; (3) \textit{virtual state}, which is a new approach to handling artificial infeasibility in SCP methods that preserves the shapes of the constraint sets; and, (4) the use of the proportional-integral projected gradient method (\pipg), a high-performance first-order conic optimization solver, in tandem with the penalized trust region (\ptr) SCP algorithm. We demonstrate the efficacy and real-time capability of {\seco} by solving a relevant multi-phase rocket landing guidance problem with nonlinear dynamics and convex constraints only, and observe that our solver is 2.7 times faster than a state-of-the-art convex optimization solver.
\end{abstract}

\end{frontmatter}


\section{INTRODUCTION}

Rocket landing guidance can be considered to be the generalization of powered-descent guidance (PDG) to include the unpowered phase(s) of flight. Precision landing techniques for orbital rockets harness convex optimization for real-time trajectory generation \citep{behcet2007jgcd, blackmore2016autonomous}. One such method, known as \textit{lossless convexification}, was the first convex optimization-based algorithm to compute a rocket landing guidance trajectory for a mid-flight large divert maneuver onboard the vehicle. Sequential convex programming (SCP) techniques for rocket landing have recently emerged as a way to handle more generalized nonconvexities in the dynamics, state constraints, environmental constraints, and path constraints (\cite{mao2016cdc}). Such SCP algorithms have been demonstrated to successfully solve a wide range of 6-DoF rocket landing problems (\cite{szmuk2019compoundSTC,reynolds2020dual,szmuk2020successive}).

Pseudospectral methods have been developed in the last decade for solving multi-phase trajectory optimization problems, particularly in aerospace applications. These methods typically parameterize the state time-history with Chebyshev and Legendre polynomials, operate on nonuniform time-grids with node points that are roots of these polynomials, and use linking conditions to tie together different phases of flight. A multi-phase Radau pseudospectral method was introduced in \citep{garrido2021ascent} for ascent and powered-descent guidance. The SPARTAN software package described in \citep{sagliano2021spartan} has been demonstrated to solve multi-phase problems arising in space applications. Further, relevant details showcasing pseudospectral approaches are provided in \citep{hwang2022integrated,ma2019direct,zhang2022pyscp}. Despite their ability to solve multi-phase problems, pseudospectral methods are typically untenable for real-time implementations \citep{malyuta2019discretization}.

In this paper, we present sequential conic optimization ({\seco}), a novel matrix-inverse-free paradigm for solving nonconvex optimal control problems in real-time, using which we formulate and solve a multi-phase rocket landing guidance problem in a unified manner. The guidance algorithm has the ability to perform mid-burn engine switching after ignition. The primary features of the {\seco} framework are as follows: (1) time-interval dilation, which allows for nonuniform time-grids and free-phase-transition-time; (2) a new formulation of compound state-triggered constraints—called \textit{single-crossing} compound state-triggered constraints—that is convex, provided the trigger and constraint conditions are convex; (3) a new approach to handling artificial infeasibility, i.e., \textit{virtual state}, that preserves the shapes of the constraint sets and does not alter the dynamics manifold; and (4) {\pipg} (proportional-integral projected gradient), a high-performance first-order solver that effectively exploits the structure of trajectory optimization problems, making it well-suited for embedded applications (\cite{yu2020proportional}). Further, we use an inverse-free exact discretization method to generate dynamically feasible solutions. Despite the layers of apparent complexity described, the method that we present is implemented via a low-footprint codebase that is easy to verify and validate.
\section{DYNAMICS}
\vspace{-0.25em}
\subsection{Time-interval dilation}
\vspace{-0.25em}
The original nonlinear dynamics, governing the evolution of state ${\f{x}{t}} \in \R^{n_{x}}$ with control input $\f{u}{t} \in \R^{n_{u}}$, over the entire time-horizon, are given by Equation \eqref{eq:nonlinear}.
\begin{align}
    \f{\dot{x}}{t} = {\f{f}{t, {\f{x}{t}}, {\f{u}{t}}}}, \enskip t \in [0, t_{f}) \label{eq:nonlinear}
\end{align}
Now, we consider the dynamics in the sub-interval $[t_{k}, t_{k+1})$, where $0 < t_{k} < t_{k+1} < t_{N} = t^{-}_{f}$, $k = \range{1}{N\!-\!1}$—where $\range{a}{b}$ denotes the range of integers between (and including) integers $a$ and $b$—and define an affine map, ${\f{\tau_{k}}{t}}$, as shown in Equation \eqref{eq:dilation}.
\begin{align}
    {\f{\tau_{k}}{t}} \defeq \dfrac{t - t_{k}}{t^{-}_{k+1} - t_{k}} \ni {\f{\tau_{k}}{t}}:[t_{k}, t_{k+1}) \to [0, 1) \label{eq:dilation}
\end{align}
This mapping is referred to as \textit{time-interval dilation}, as it normalizes the wall-clock time-interval to a known fixed interval—in our case $[0, 1)$—by either shrinking or expanding—and hence \textit{dilating}—the original time-interval. Next, we apply the derivative operator with respect to the dilated time $\tau_{k}$, denoted by $\circdot{\square}$, to Equation \eqref{eq:nonlinear}, and invoke the chain-rule, as shown in Equation \eqref{eq:chain}, where $t \in {[t_{k}, t_{k+1})}$.
\begin{flalign}
{\f{\circdot{x}}{t}} &= \frac{d}{d\tau\mathrlap{\!\!\:_{k}}}\,
{\f{x}{t}} = \frac{d t}{d\tau\mathrlap{\!\!\:_{k}}}~\frac{d}{d t}\,
{\f{x}{t}} = \frac{d t}{d\tau\mathrlap{\!\!\:_{k}}}~{\f{\dot{x}}{t}} = \underbracket[0.14ex]{\left(t^{-}_{k+1}\! - t_{k}\right)}_{s_{k}}\,{\f{\dot{x}}{t}} \nonumber\\
&= s_{k}\,{\f{f}{t, {\f{x}{t}}, {\f{u}{t}}}} \defeq {\f{F}{t, {\f{x}{t}}, {\f{u}{t}}, s_{k}}} \label{eq:chain}
\end{flalign}
The multiplier in Equation \eqref{eq:chain}, $s_{k} \defeq t^{-}_{k+1} - t_{k} \in \R_{+}$, which is nothing but the length of the $k\textsuperscript{th}$ wall-clock time-interval, is termed the \textit{dilation factor}. By treating a phase-based subset of $s_{k},\,k = \range{1}{N\!-\!1}$, i.e., such that $k \subseteq \range{1}{N\!-\!1}$, as decision variables and discretizing the system over them, we allow the optimizer to decide what the temporal spacing of discrete nodes should be in each phase rather than use a uniform temporal grid over the entire horizon. In the approach we propose, this is the key to enabling free-transition-time multi-phase trajectory optimization within a single-shot optimization framework, without requiring any mixed-integer or nonconvex constraints to handle the discrete temporal events/switching.

Although it is possible to allow each dilation factor to be an independent decision variable, we choose to partition the temporal grid based on the phases of flight, and evenly space the temporal nodes within each phase. This measure is taken to mitigate extreme inter-sample constraint violation, which tends to occur when fully adaptive grids are used. The time-dilated dynamics given by Equation \eqref{eq:chain} will be used henceforth. Note that using the time-dilated dynamics given by Equation \eqref{eq:chain} in lieu of Equation \eqref{eq:nonlinear} converts the original \textit{free-final-time} optimal control problem to an equivalent \textit{fixed-final-time} optimal control problem, with the effective horizon being $[0, N\!-\!1)$.

\vspace{-0.25em}
\subsection{Linearization}
\vspace{-0.25em}

A convex approximation of the original nonconvex optimal control problem is obtained by linearizing the nonlinear dynamics (Equation \eqref{eq:chain})—which leads to a linear time-varying (LTV) system—and keeping the convex constraints intact. The state, control, and parameter constraint sets are assumed to belong to a set that has separable closed-form projection operations. We stress that most of the common constraints in trajectory optimization problems naturally fit this template \citep{SCPToolboxCSM2022}. Nonconvex constraints can be handled using this framework too, by means of either linearization or convex approximation.

\vspace{-0.25em}
\subsection{Discretization}
\vspace{-0.25em}

We assume a first-order hold (FOH) on the control input signal throughout the horizon and a zero-order hold (ZOH) on the control input signal only during discrete switching events such as engine startup and downselection. Both of these control parameterizations have two key properties that make them attractive for optimal control applications: (1) inter-sample satisfaction of the convex control constraints is guaranteed (provided they are satisfied at the discrete temporal nodes), which is in contrast to pseudospectral methods; and, (2) the resulting conic subproblem has a sparsity pattern that is amenable to real-time implementation \citep{malyuta2019discretization, szmuk2020successive}.

In the FOH case, the control profile is parameterized as shown in Equation \eqref{eq:control_basis_func}, where $t \in [t_{k}, t_{k+1})$ and $\tau_{k} \in [0, 1)$, as given by Equation \eqref{eq:dilation}.
\begin{gather}
    \f{u}{\tau_{k}} = {\left(1 - \tau_{k}\right)}\,u_{k} + \tau_{k}\,u_{k+1}, \enskip k = \range{1}{N\!-\!1}
    \label{eq:control_basis_func}
\end{gather}
The LTV dynamics can now be written: (1) using the piecewise-affine control input parameterization given by Equation \eqref{eq:control_basis_func}; and, (2) in terms of deviations from the reference, as shown in Equation \eqref{eq:dyn_lin_interp}. The reference quantities are denoted by $\overline{\square}$, and $\Delta \square$ denotes the deviation of a quantity from its reference, i.e., $\Delta \square \defeq \square - \overline{\square}$, and $\f{\Delta \circdot{x}}{\tau_{k}} \defeq \f{\circdot{x}}{\tau_{k}} - \f{F}{\tau_{k}, {\f{\overline{x}}{\tau_{k}}}, {\f{\overline{u}}{\tau_{k}}}, \overline{s}_{k}}$. The approximate nature of the equation is an artifact of linearization of the original nonlinear dynamics via truncation of the higher-order ($\ge 2$) terms in the Taylor series expansion.
\begin{align}
    \begin{split}
    {\f{\Delta\accentset{\circ}{x}}{\tau_{k}}} &\approx {\f{A}{\tau_{k}}}\,{\f{\Delta x}{\tau_{k}}} + {\f{B}{\tau_{k}}}\, {\left(1 - \tau_{k}\right)}\,\Delta u_{k} \\
    &\hphantom{\approx\;}+ {\f{B}{\tau_{k}}}\,\tau_{k}\,\Delta u_{k+1} + {\f{S}{\tau_{k}}}\,\Delta s_{k}
    \end{split}
    \label{eq:dyn_lin_interp}
\end{align}
The \textit{state transition matrix} (STM) associated with Equation \eqref{eq:dyn_lin_interp}, denoted by ${\f{\Phi}{\tau_{k}, 0}}$, $\tau_{k} \in [0, 1)$, satisfies the following matrix differential equation: \({\f{\circdot{\Phi}}{\tau_{k}, 0}} = {\f{A}{\tau_{k}}}\,{\f{\Phi}{\tau_{k}, 0}}\), with ${\f{\Phi}{0, 0}} = I_{n_{x}}$. The unique solution to Equation \eqref{eq:dyn_lin_interp} is given by Equation \eqref{eq:lin_sys_theory} \citep{antsaklis2006linear, SCPToolboxCSM2022}.
\begin{alignat}{1}
        &{\f{\Delta x}{\tau_{k}}} = {\f{\Phi}{\tau_{k}, 0}}\,{\f{\Delta x}{0}} + \int_{0}^{\tau_{k}} {\f{\Phi}{\tau_{k}, \zeta}}\,\cdot \label{eq:lin_sys_theory}\\
        & \cdot{\left\{{\f{B}{\zeta}}\,{\left(1 - \tau_{k}\right)}\,\Delta u_{k} + {\f{B}{\zeta}}\,\tau_{k}\,\Delta u_{k+1} + {\f{S}{\zeta}}\,\Delta s_{k}\right\}}\,\mathrm{d}\zeta \nonumber
\end{alignat}
\begin{figure}[H]
    \centering
    \vspace{-0.25em}
    \includegraphics[width=0.88125\linewidth]{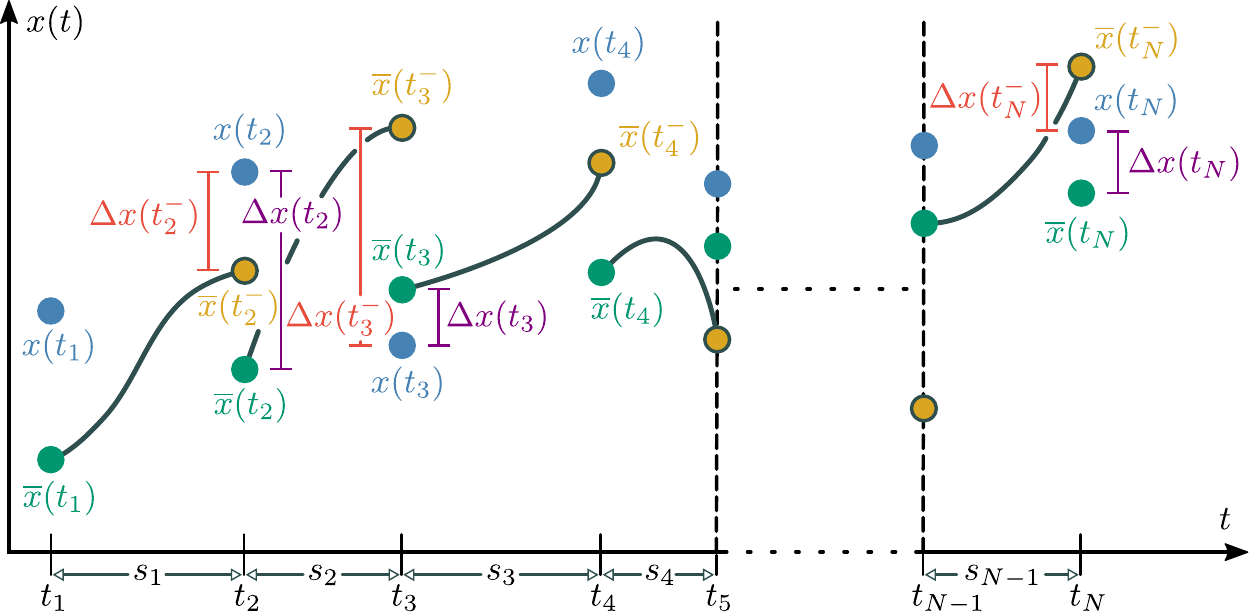}
    \vspace{-0.5em}
    \caption{Propagation of the state. ${\f{\Delta x}{t_{1}}}\!\defeq 0$, where $t_{1}\!\defeq 0$. The \textit{stitching condition} for $k = \range{1}{N\!-\!1}$ is given by ${\color{bricks}{\f{\Delta x}{t^{-}_{k+1}}}} + {\color{goldenrod}{\f{\overline{x}}{t^{-}_{k+1}}}} = {\color{purple}{\f{\Delta x}{t_{k+1}}}} + {\color{darts}{\f{\overline{x}}{t_{k+1}}}} = {\color{steelblue}{\f{x}{t_{k+1}}}}$.}
    \label{fig:stitching}
    \vspace{-0.75em}
\end{figure}
Evaluating Equation \eqref{eq:lin_sys_theory} at $\tau_{k} = 1^{-}$, we get Equation \eqref{eq:dyn_disc}. We replace the limits $0$ and $1$ with $0_{k}$ and $1_{k}$, respectively, to explicitly indicate the dependence on index $k$, i.e., ${\f{\tau_{k}}{t_{k}}} \defeq 0_{k}$ and ${\f{\tau_{k}}{t^{-}_{k+1}}} \defeq 1^{-}_{k}$.
\begin{flalign}
    \hspace{-0.675em}\f{\Delta x}{1_{k}^{-}} &= A_{k} \f{\Delta x}{0_{k}} + B_{k}^{-} \Delta u_{k} + B_{k}^{+} \Delta u_{k+1} + S_{k} \Delta s_{k}
\label{eq:dyn_disc}%
\end{flalign}
$A_{k}$, $B_{k}^{-}$, $B_{k}^{+}$, and $S_{k}$ can be computed as the solution to the initial value problem (IVP) given by Equations \eqref{eqs:disc_ivp}, respectively, integrated from $0_{k}$ to $1_{k}^{-}$.
\begin{subequations}
\begin{align}
    {\f{\circdot{\Psi}_{A}}{\zeta}} &= {\f{A}{\zeta}}\,{\f{\Psi_{A}}{\zeta}} \\
    {\f{\circdot{\Psi}_{B^{-}}}{\zeta}} &= {\f{A}{\zeta}}\,{\f{\Psi_{B^{-}}}{\zeta}} + {\f{B}{\zeta}}\,{\left(1 - \zeta\right)} \label{eq:Bkm}\\
    {\f{\circdot{\Psi}_{B^{+}}}{\zeta}} &= {\f{A}{\zeta}}\,{\f{\Psi_{B^{+}}}{\zeta}} + {\f{B}{\zeta}}\,\zeta \\
    {\f{\circdot{\Psi}_{S}}{\zeta}} &= {\f{A}{\zeta}}\,{\f{\Psi_{S}}{\zeta}} + {\f{S}{\zeta}}
\end{align}
\label{eqs:disc_ivp}%
\end{subequations}
where function \({\f{\Psi_{A}}{\zeta}}\) is defined such that \( \zeta \mapsto {\f{\Phi}{\zeta, 0_{k}}}\), and the initial conditions for Equations \eqref{eqs:disc_ivp} are: ${\f{\Psi_{A}}{0_{k}}} = I_{n_{x}}$, ${\f{\Psi_{B^{-}}}{0_{k}}} = {\f{\Psi_{B^{+}}}{0_{k}}} = 0_{n_{x} \times n_{u}}$, and ${\f{\Psi_{S}}{0_{k}}} = 0_{n_{x}}$.
Note that Equation \eqref{eq:stitching}, which we refer to as the \textit{stitching condition}, holds, as is evident from Figure \ref{fig:stitching}.
\begin{align}
    \f{\Delta x}{1_{k}^{-}} + \f{\overline{x}}{1_{k}^{-}} = \f{\Delta x}{1_{k}} + \f{\overline{x}}{1_{k}} \label{eq:stitching}
\end{align}
The discretized dynamics can now be given by Equation \eqref{eq:dyn_disc_stitch}, where $\Delta x_{k} \defeq \f{\Delta x}{0_{k}}$, $x_{k+1}^{\mathrm{prop}} \defeq \f{\overline{x}}{1^{-}_{k}}$, $\overline{x}_{k+1} \defeq \f{\overline{x}}{1_{k}}$, and $\f{\overline{u}}{\tau_{k}} \defeq {\left(1 - \tau_{k}\right)}\,\overline{u}_{k} + \tau_{k}\,\overline{u}_{k+1}$, for $1 \le k \le N\!-\!1$.
\begin{align}
\begin{split}
    \Delta x_{k+1} = A_{k} \Delta x_{k} + B_{k}^{-} \Delta u_{k} &+ B_{k}^{+} \Delta u_{k+1} + S_{k} \Delta s_{k} \\
    &+ x_{k+1}^{\mathrm{prop}} - \overline{x}_{k+1}
\end{split}
\label{eq:dyn_disc_stitch}%
\end{align}
The discretized dynamics in terms of the absolute variables are recovered from Equation \eqref{eq:dyn_disc_stitch}, as shown in Equation \eqref{eq:dyn_disc_abs}.
\begin{align}
    x_{k+1} &= A_{k}\,x_{k} + B_{k}^{-}\,u_{k} + B_{k}^{+}\,u_{k+1} + S_{k}\,s_{k}\; + \label{eq:dyn_disc_abs}\\ 
    &\hphantom{=~} x_{k+1}^{\mathrm{prop}} - \left(A_{k}\,\overline{x}_{k} + B_{k}^{-}\,\overline{u}_{k} + B_{k}^{+}\,\overline{u}_{k+1} + S_{k}\,\overline{s}_{k}\right) \nonumber
\end{align}
Equation \eqref{eq:dyn_disc_abs} represents an \textit{exact} discretization of the LTV dynamics, which means that the error between the continuous-time trajectory and the discrete-time trajectory at the discrete temporal nodes is analytically zero.

In order to obtain the corresponding expressions for the ZOH case, the following changes are incorporated: (1) Equation \eqref{eq:Bkm} is replaced by Equation \eqref{eq:Bk} (and accordingly, $B_{k}^{-}$ is replaced by $B_{k}$); and, (2) $B_{k}^{+}$ is set to a zero matrix with the same dimensions.
\begin{align}
    {\f{\circdot{\Psi}_{B}}{\zeta}} &= {\f{A}{\zeta}}\,{\f{\Psi_{B}}{\zeta}} + {\f{B}{\zeta}} \label{eq:Bk}
\end{align}
\section{STATE-TRIGGERED CONSTRAINTS}
\vspace{-0.25em}

We introduce a specialized formulation of compound state-triggered constraints (STCs) for problems in which the trigger functions are activated only once and are strictly monotonic in a neighborhood around which they are activated. We refer to these constraints as \textit{single-crossing} compound state-triggered constraints. 

Let $g(\cdot)$ be a trigger function that is said to be \textit{activated} on the set $\{x\,|\,g(x)\le g^\star\}$ for some trigger value $g^\star$. Further, let $g^{-1}(g^\star)$ be a well-defined pre-image. Then, $g(\cdot)$ is called a \textit{single-crossing} trigger function if $x^\star \defeq g^{-1}(g^\star)$ is a singleton and $g(\cdot)$ is strictly monotonic in
a neighborhood around $x^\star$. Such a formulation is especially useful in applications such as rocket landing that require certain STCs to be satisfied for mission success, wherein it is reasonable to expect the trigger conditions to be activated once and only once. For instance, it is reasonable to expect/require a rocket in descent from a certain initial altitude, $h_{i}$, with its target landing location at the origin, to certainly cross a trigger altitude $0 < h_{\text{trigger}} < h_{i}$ once during its descent, and not surpass that altitude after.

For the purpose of demonstration, we consider a compound STC with one trigger condition and multiple constraint conditions to be imposed with the \textsc{and} logic. Let $\f{g}{\cdot}$ be the trigger function and $\f{c^{j}}{\cdot}$, $j = \range{1}{n_{c}}$, be the constraint functions, where $n_{c}$ is the number of constraint conditions. The purpose of the compound STC \cite[Section II.B]{szmuk2019compoundSTC} is to satisfy the condition given in Equation \eqref{eq:STC_condition}, where $\f{x}{t}$ is the state. $\forall t \in [0, t_{f})$,
\begin{align}
\f{g}{\f{x}{t}} \le 0 \Rightarrow \bigwedge_{j=1}^{n_c} \f{c^{j}}{\f{x}{t}} \le 0
    \label{eq:STC_condition}
\end{align}
In maneuvers with a single-crossing trigger condition, the trigger condition is activated once and only once, i.e., $\f{g}{\f{x}{t}} = 0$ is guaranteed to activate at some $t = t_{\mathrm{trigger}}$, $\f{g}{\f{x}{t}} > 0$ $\forall t \in [0, t_{\mathrm{trigger}})$, and $\f{g}{\f{x}{t}} < 0$ $\forall t \in (t_{\mathrm{trigger}}, t_{f})$. Using this fact, the single-crossing compound STC is formulated as shown in Equation \eqref{eq:SCCSTC}. We emphasize that $t_\mathrm{trigger}$ itself is free, and hence the compound constraint is state-triggered and not time-triggered.
\begin{subequations}
\begin{align}
    \begin{aligned}\f{g}{\f{x}{t}} \ge 0 \end{aligned}\quad\enskip\,\:
    &\forall t \in [0, t_{\mathrm{trigger}}) \\
    \begin{aligned}\f{g}{\f{x}{t}} = 0 \end{aligned},\quad\enskip\!\;
    &\hphantom{\forall} t = t_{\mathrm{trigger}} \label{eq:waypoint} \\
    \left.\begin{aligned}\f{g}{\f{x}{t}} &\le 0 \\
     \f{c^{1}}{\f{x}{t}} &\le 0 \\
     &\shortvdotswithin{=}
     \f{c^{n_{c}}}{\f{x}{t}} &\le 0 \end{aligned}\enskip
     \right\rbrace \enskip &\forall t \in (t_{\mathrm{trigger}}, t_{f})  
\end{align}\label{eq:SCCSTC}%
\end{subequations}
If the trigger conditions and the constraint conditions above are individually convex, the compound STC is entirely convex. This is in contrast to existing formulations of STCs in the literature that are inherently nonconvex, regardless of the convexity of the trigger and constraint conditions \citep{szmuk2019compoundSTC, szmuk2020successive, reynolds2020dual}. However, those methods are more general, in that they do not mandate the trigger condition to be single-crossing—we trade off this generality for convexity/simplicity in our method. Further, we note that the single-crossing STC formulation bears resemblance to the method adopted in \citep{bhasin2016fuel}.

The aforementioned approaches in the literature use uniform temporal spacing of discrete nodes over the entire horizon. As a result, along with the fact that inter-sample constraint satisfaction is typically not guaranteed in general, these approaches do not guarantee the imposition of the constraints exactly at the specified trigger conditions. The triggering of these constraints (in time) is only accurate up to the spacing of the grid, and the corresponding solutions usually demonstrate violation of these constraints with respect to the triggers. This issue becomes more prevalent in maneuvers over very long time horizons, especially when the set of feasible trigger windows is much smaller in comparison, and can be detrimental to mission success if accurate triggering is required. 

We leverage time-interval dilation to impose single-crossing STCs, and treat the windows within which these constraints need to be imposed as distinct phases of flight. This allows for a fine grid in phases that involve critical constraints that need to be satisfied to ensure mission success, and a coarser grid in the more benign phases of flight—thus enabling one-shot multi-phase trajectory optimization. If the solution converges to a feasible trajectory, the STCs are guaranteed to be satisfied at the triggers (since the triggers are imposed as waypoints as in Equation \eqref{eq:waypoint}), and at every discrete temporal node within the trigger window.
\section{SEQUENTIAL CONIC OPTIMIZATION}

\subsection{Virtual state}\label{subsec:virtual_state}

\textit{Artificial infeasibility} refers to the phenomenon wherein a subproblem can become infeasible as a result of linearization of the nonconvex constraints even if there exists a feasible solution to the original problem. Typically, unconstrained yet heavily penalized slack variables are added to the linearized constraints, so as to ensure that the subproblem is always feasible.

We propose a new approach to handling artificial infeasibility by means of a \textit{virtual state} variable, which serves as a copy of the original state. This approach helps decouple the dynamics and control constraints from the state constraints and exactly satisfy all the path constraints at each solver iteration, while ensuring that the subproblem never turns infeasible. If $x$ is the actual state variable, $u$ is the control variable, and $\xi$ is the virtual state variable, the dynamics constraint is imposed on $x$ and $u$, the control constraints are imposed on $u$, and all the state constraints are imposed on $\xi$.

To ensure that the dynamics and all other constraints are satisfied at convergence, we minimize the error between $x$ and $\xi$ by heavily penalizing the squared distance between them in the objective function. The virtual state does not alter the dynamics manifold (unlike the virtual control approach \citep{szmuk2020successive, reynolds2020dual}), and preserves the shapes of the conic state constraint sets (unlike the virtual buffer approach \citep{SCPToolboxCSM2022}).
\vspace{-0.125em}

\subsection{Conic subproblem}

We impose a soft trust region on the decision variable and use the penalized trust region (\ptr) algorithm \citep{szmuk2020successive, reynolds2020real, reynolds2020dual}. The discretized conic subproblem with virtual state(s) and a soft trust region is shown in Problem \eqref{prob:conic_subproblem}, which is strongly convex.
\begin{subequations}
\begin{alignat}{4}
    \underset{u,\,s}{\text{min}} \quad &w_{\mathrm{c}} {\f{J}{x_{N}}} + \tfrac{1}{2}\!\left(w_{\mathrm{tr}} J_{\mathrm{tr}} + w_{\mathrm{vse}} J_{\mathrm{vse}}\right) \span\span\span\span\label{eq:objective}\\
        \text{s.t.} \quad
        & x_{k+1} = \text{RHS of Eq. \eqref{eq:dyn_disc_abs}},\span \enskip &k = \range{1}{N\!-\!1}\label{eq:dyn_constr_disc}\\
        & \xi_{k} \in \mathcal{X}_{\sector}, \span\enskip &k = \range{1}{N} \label{eq:state_constr_disc}\\
        & u_{k} \in \mathcal{U}_{\sector}, \span\enskip &k = \range{1}{N} \label{eq:control_constr_disc}\\
        & s_{k} \in \mathcal{S}_{\sector}, \span\enskip &k = \range{1}{n_\mathrm{phase}}\label{eq:param_constr_disc}
    \end{alignat}
\label{prob:conic_subproblem}%
\end{subequations}
where $x \in \R^{n_{x}}$ is the state vector, $\xi \in \R^{n_{x}}$ is the \textit{virtual state} vector, $u \in \R^{n_{u}}$ is the control input vector, and $s \in \R_{++}$ is the dilation factor (vector); $\mathcal{X}_{\sector}$, $\mathcal{U}_{\sector}$, and $\mathcal{S}_{\sector}$ are the state, control, and temporal constraint sets, respectively, which are assumed to be closed and convex; ${\f{J}{x_{N}}}$ is the original cost function, assumed to be in the Mayer form \citep{berkovitz2013optimal},
$J_{\mathrm{tr}} \defeq \sum_{k=1}^{N} \!\left(\norm{x_{k} - \overline{x}_{k}}_{2}^{2} + \norm{u_{k} - \overline{u}_{k}}_{2}^{2}\right) + \sum_{k=1}^{n_{\mathrm{phase}}} \norm{s_{k} - \overline{s}_{k}}_{2}^{2}$ is the trust region penalty, and $J_{\mathrm{vse}} \defeq \sum_{k=1}^{N} \norm{x_{k} - \xi_{k}}_{2}^{2}$ is the virtual state error penalty, $w_{c}$, $w_{\mathrm{tr}}$, and $w_{\mathrm{vse}}$ being their respective weights. 
\vspace{0.375em}

\begin{prop}
The virtual state error penalty, given by $\sum_{k=1}^{N} \norm{x_{k} - \xi_{k}}_{2}^{2}$, is convex.
\end{prop}
\begin{pf}
Let $y_{k} \defeq \left(x_{k},\,\xi_{k}\right)$. \,$\therefore$ $\norm{x_{k} - \xi_{k}}_{2}^{2} = y_{k}^{\top} M y_{k}$, where $M = \begin{pmatrix}\hphantom{-}1\hphantom{.} & -1 \\ -1\hphantom{.} & \hphantom{-}1\end{pmatrix} \otimes I_{n_{x}}$. Since $\operatorname{spec} M \in \{0, 2\}$, $M$ is positive semidefinite (PSD). Therefore, the quadratic form, $y_{k}^{\top} M y_{k}$, is convex. Since the sum of convex functions is convex, $\sum_{k=1}^{N} y_{k}^{\top} M y_{k} = \sum_{k=1}^{N} \norm{x_{k} - \xi_{k}}_{2}^{2}$ is also convex. \hfill $\blacksquare$
\end{pf}

\vspace{0.25em}

Problem \eqref{prob:conic_subproblem} can be vectorized, i.e., assembled into the form of Problem \eqref{prob:conic}, by stacking all the decision variables into a single vector, $z$. For more details, see \citep{yu2022extrapolated}.
\begin{subequations}
\begin{align}
    \underset{z}{\text{min}} \quad &\frac{1}{2} z^{\top} Q z + \langle q, z \rangle \label{eq:vec_obj}\\
    \text{s.t.} \quad
        & H z - h = 0 \label{eq:vec_dyn}\\
        & z \in \D \label{eq:vec_proj}
\end{align}
\label{prob:conic}%
\end{subequations}
\vspace{-1em}
\section{HIGH-PERFORMANCE SOLVER}
\vspace{-0.375em}

The proportional-integral projected gradient method (\pipg) \citep{yu2020proportional,yu2021proportional} is a first-order primal-dual optimization algorithm for solving conic optimization problems such as the one given by Equations \eqref{prob:conic}. {\pipg} is compatible with extrapolation, which has been shown to improve its practical convergence behavior. Here, we use the extrapolated {\pipg} algorithm, also denoted by {\begin{footnotesize}x\end{footnotesize}\pipg}, given by the iterative sequence in Equation \eqref{alg:pipg} \citep{yu2022extrapolated}.
\begin{subequations}
\begin{align}
z^{j+1} &= \pi_{\mathbb{D}}\!\left[\zeta^j-\alpha\!\left(Q\zeta^j+q+H^{\top}\!\eta^j\right)\right] \\
w^{j+1} &= \eta^j+\beta\!\left(H\!\left(2 z^{j+1}-\zeta^j\right)-h\right) \\
\zeta^{j+1} &= (1-\rho) \zeta^j+\rho z^{j+1} \\
\eta^{j+1} &= (1-\rho) \eta^j+\rho w^{j+1}
\end{align}\label{alg:pipg}%
\end{subequations}
where the primal variable, $z^j$, converges to an optimum of Problem \eqref{prob:conic} as $j\to\infty$. The step-sizes, $\alpha$ and $\beta$, are computed in accordance with \cite[Lemma 2]{yu2022extrapolated}, and $\rho \in [1, 2)$ is the extrapolation factor. We leverage the warm-starting capability of {\pipg} and use the $\ell_{2}$-hypersphere preconditioning technique described in \citep{kamath2023customized} to further accelerate convergence.
\section{MULTI-PHASE ROCKET LANDING GUIDANCE}
\vspace{-0.375em}

As a representative example for multi-phase rocket landing guidance via {\seco}, we consider the problem of terrestrial precision landing of a vehicle akin to Starship, which is designed to be a fully reusable rocket and currently in development \citep{shotwellspace}.

We consider a nonlinear planar model of the vehicle, given by Equations \eqref{eq:starship_model}, with $\hat{x}$ being the vertical axis and $\hat{z}$ being the horizontal axis in the inertial frame (and the longitudinal and lateral axes in the body frame, respectively). We use a simple, tractable aerodynamic model to account for the aerodynamic effects on the vehicle during descent \citep{szmuk2020successive}. In practice, however, numerical databases can be used along with {\seco} \citep{mceowen2022hypersonic}. Further, with our approach, aerodynamic-free polynomial coast-phase ballistic trajectory predictions, such as the one provided in \citep{szmuk2020successive}, are not required. We note that \citep{lee2022multi} provide an SCP-based 6-DoF multi-phase rocket landing guidance implementation, but impose nonconvex STCs and implement inexact (trapezoidal) discretization, unlike our approach.

The \textit{original} state vector is defined as follows: ${\f{x}{t}} \defeq\\ \left({\f{m}{t}}, {\f{r}{t}}, {\f{v}{t}}, {\f{\theta}{t}}, {\f{\omega}{t}}\right)$, where $\f{m}{t} \in \R_{+}$ is the mass, ${\f{r}{t}} \in \R^{2}$ is the position, ${\f{v}{t}} \in \R^{2}$ is the velocity, ${\f{\theta}{t}} \in \R$ is the body tilt angle with respect to the inertial vertical, and ${\f{\omega}{t}} \in \R$ is the angular velocity of the body.

The control input vector is defined as follows: ${\f{u}{t}} \defeq \left({\f{T}{t}}, {\f{\delta}{t}}\right)$, where ${\f{T}{t}} \in \R$ is the thrust magnitude and ${\f{\delta}{t}} \in \R$ is the gimbal deflection angle.

Further, we define a \textit{virtual state} vector and impose all the state constraints, including the boundary conditions, on this variable, in accordance with Subsection \ref{subsec:virtual_state}: ${\f{\xi}{t}} \defeq \left({\f{m^{\xi}}{t}}, {\f{r^{\xi}}{t}}, {\f{v^{\xi}}{t}}, {\f{\theta^{\xi}}{t}}, {\f{\omega^{\xi}}{t}}\right)$.
\begin{subequations}
\begin{align}
    {\f{\dot{m}}{t}} &= -\alpha_{\mathrm{e}}\!\:{\f{T}{t}} \\
    {\f{\dot{r}}{t}} &= {\f{v}{t}} \\
    {\f{\dot{v}}{t}} &= \frac{1}{\f{m}{t}}\!\left({\f{F_{\I}}{t}} + {\f{A_{\I}}{t}}\right) + g \\
    {\f{\dot{\theta}}{t}} &= {\f{\omega}{t}} \\
    {\f{\dot{\omega}}{t}} &= \frac{1}{\f{J}{t}}\!\left({\f{F_{\B_{\hat{z}}}\!}{t}}\,l_{\mathrm{cm}} - {\f{A_{\B_{\hat{z}}}\!}{t}}\,l^{j}_{\mathrm{cp}} \right)
\end{align}
\label{eq:starship_model}
\end{subequations}
where
\begin{subequations}
\begin{align}
    \alpha_{\mathrm{e}} &\defeq \frac{1}{I_{\mathrm{sp}}\!\:g_{0}}\\
    {\f{F_{\I}}{t}} &\defeq {\f{T}{t}}\!\begin{pmatrix}\hphantom{-\!\!\:}{\f{\cos}{\f{\theta}{t} + {\f{\delta}{t}}}} \\ -{\f{\sin}{\f{\theta}{t} + {\f{\delta}{t}}}} \end{pmatrix} \\
    {\f{A_{\I}}{t}} &\defeq {\f{R_{\I\leftarrow\B}}{t}}\!\:{\f{A_{\B}}{t}} \\
    {\f{F_{\B}}{t}} &\defeq {\f{T}{t}}\!\begin{pmatrix}\hphantom{-\!\!\:}{\f{\cos}{\f{\delta}{t}}} \\ -{\f{\sin}{\f{\delta}{t}}} \end{pmatrix} \\
    {\f{A_{\B}}{t}} &\defeq -\frac{\rho_{\mathrm{air}}\!\:S_{\mathrm{area}}\!\:\norm{\f{v}{t}}_{2}}{2} C_{\mathrm{aero}}\!\:{\f{R_{\I\leftarrow\B}^{\top}}{t}}\!\:{\f{v}{t}} \\
    {\f{R_{\I\leftarrow\B}}{t}} &\defeq \begin{pmatrix}\hphantom{-\!\!\:}\cos{\f{\theta}{t}} & \hphantom{\!\!\:}\sin{\f{\theta}{t}}\\-\sin{\f{\theta}{t}} & \cos{\f{\theta}{t}}\end{pmatrix}
\end{align}
\end{subequations}
Here, $\alpha_{\mathrm{e}} \in \R_{++}$ is the thrust-specific fuel consumption (TSFC), $I_{\mathrm{sp}} \in \R_{++}$ is the specific impulse of the rocket engine, $g_{0} \in \R_{++}$ is standard Earth gravitational acceleration, $g \defeq (-g_{0}, 0)$, ${\f{R_{\I\leftarrow\B}}{t}} \in {\f{\mathrm{SO}}{2}}$ is the rotation matrix that maps coordinates in the body frame to the inertial frame, $\rho_{\mathrm{air}} \in \R_{++}$ is the ambient atmospheric density, $S_{\mathrm{area}} \in \R_{++}$ is the reference area, $C_{\mathrm{aero}} \defeq \operatorname{diag}\{c_{\hat{x}}, c_{\hat{z}}\}$ is the aerodynamic coefficient matrix, where $c_{\hat{x}}, c_{\hat{z}} \in \R_{++}$ are the aerodynamic coefficients along the body $\hat{x}$ and $\hat{z}$ axes, respectively, and ${\f{J}{t}} \in \R_{++}$ is the moment of inertia of the vehicle about the body $\hat{y}$ axis (out-of-plane). The body of the vehicle is assumed to be a uniform solid cylinder, and hence, the moment of inertia about its central diameter is given by ${\f{J}{t}} \defeq {\f{m}{t}}\!\left(\frac{l_{\mathrm{r}}^{2}}{4} + \frac{l_{\mathrm{h}}^{2}}{12}\right)$, where $l_{\mathrm{r}} \in \R_{++}$ and $l_{\mathrm{h}} \in \R_{++}$ are the radius and height of the fuselage, respectively.

The engines are assumed to be co-located, and the location of the vehicle mass-center is assumed to be fixed in the body frame. The thrust moment-arm (the distance between the vehicle mass-center and the engine gimbal hinge point) is denoted by $l_{\mathrm{cm}} \in \R_{++}$, and $l^{j}_{\mathrm{cp}} \in \R_{+}$ is the aerodynamic moment-arm (the distance between the vehicle mass-center and the center-of-pressure), where $j \in \{0, 1\}$; $j = 0$ for the unpowered phase of flight and $j = 1$ for the powered phases of flight—$l^{0}_{\mathrm{cp}}$ is assumed to be maintained at zero via independent aerodynamic controls by means of forward and aft flaps (aerodynamic control surfaces), i.e., the center-of-pressure and the mass-center are assumed to be coincident when the vehicle is in the coast phase.

At the engine ignition (PDI) epoch, it is assumed that the following events occur: (1) the forward flaps are fully extended (to maximize drag towards the nose-cone); and, (2) the aft flaps are fully folded (to minimize drag towards the aft section of the vehicle). As a result, the center-of-pressure shifts away from the mass-center, towards the nose-cone, and induces an aerodynamic torque (pitching moment) on the vehicle. Hence, $l^{1}_{\mathrm{cp}}$ is set to a nonzero value, and it is kept fixed for the remainder of the trajectory. This, along with gimbaling of the rocket engines, is used to induce the flip maneuver to get the vehicle upright in preparation for terminal descent.

The guidance problem is partitioned into four phases: (1) the unpowered, subsonic coast phase; (2) the high-thrust (3-engine) burn phase; (3) the low-thrust (1-engine) burn phase; and, (4) the altitude-triggered terminal descent phase, separated by the following important discrete events/epochs: (1) powered-descent initiation (PDI) or engine ignition, $t_{\mathrm{ignition}}$; (2) engine downselection or switching from a triple-engine burn to a single-engine burn, $t_{\mathrm{switch}}$; and, (3) altitude-based triggering of the terminal descent phase, $t_{\mathrm{trigger}}$. We emphasize that putting all of these phases together implicitly leads to a free-ignition-time, free-engine-switching-time, \textit{and} free-final-time optimal control problem (subject to the temporal constraints imposed).

For the discretized problem, the temporal grid we choose is as given by Equation \eqref{eq:discrete_grid}, where $N$ is the number of discrete temporal nodes and $k_{\mathrm{ignition}}$, $k_{\mathrm{switch}}$, and $k_{\mathrm{trigger}}$ are the nodes at which the discrete events occur.
\begin{align}
k \in \{1, \ldots, k_{\mathrm{ignition}}, \ldots, k_{\mathrm{switch}}, \ldots, k_{\mathrm{trigger}}, \ldots, N\}
\label{eq:discrete_grid}
\end{align}
\vspace{-1.125em}

\begin{figure*}[htp!]
    \centering
    \begin{minipage}[b]{0.25625\linewidth}
        \centering
        \includegraphics[width=1.0125\linewidth]{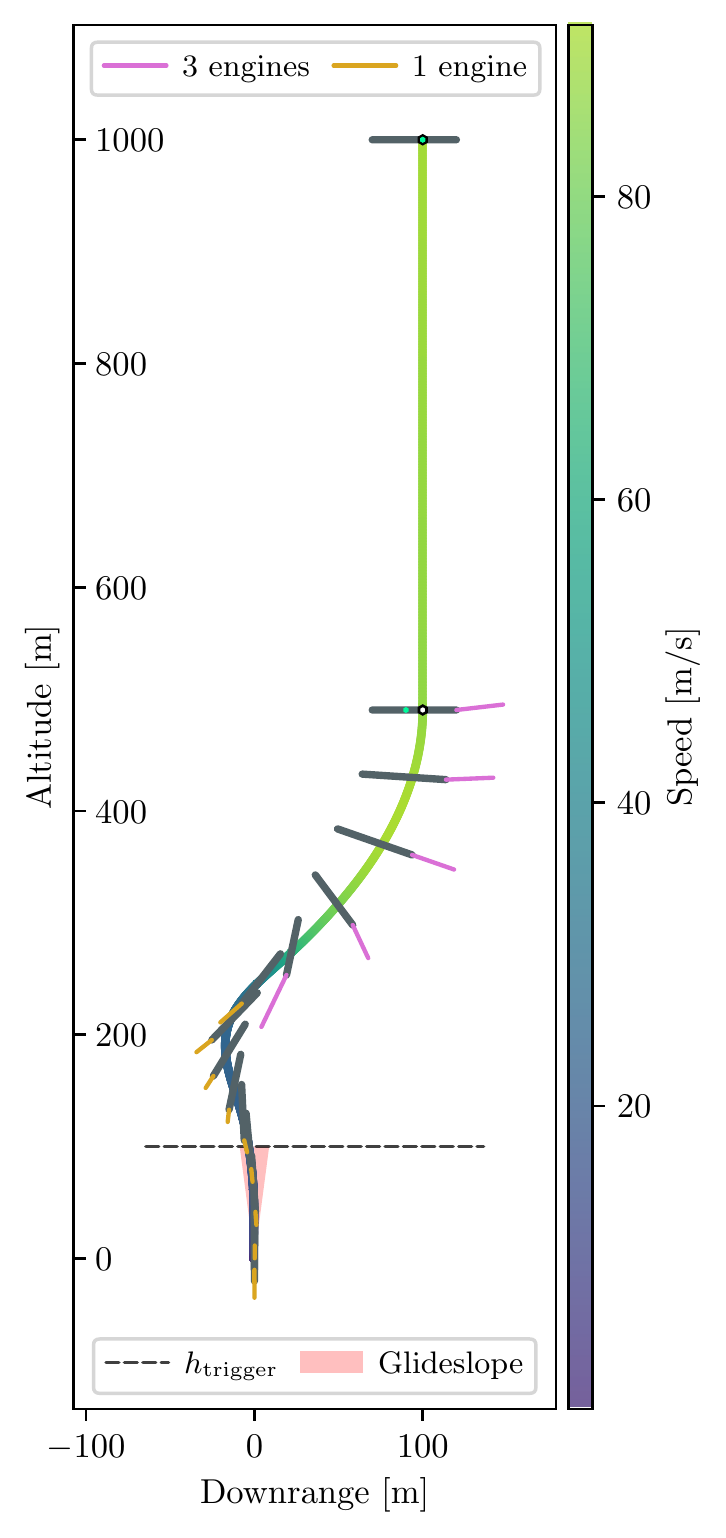}
    \end{minipage}%
    \begin{minipage}[b]{0.74375\linewidth}
        \centering
        \includegraphics[width=1.0125\linewidth]{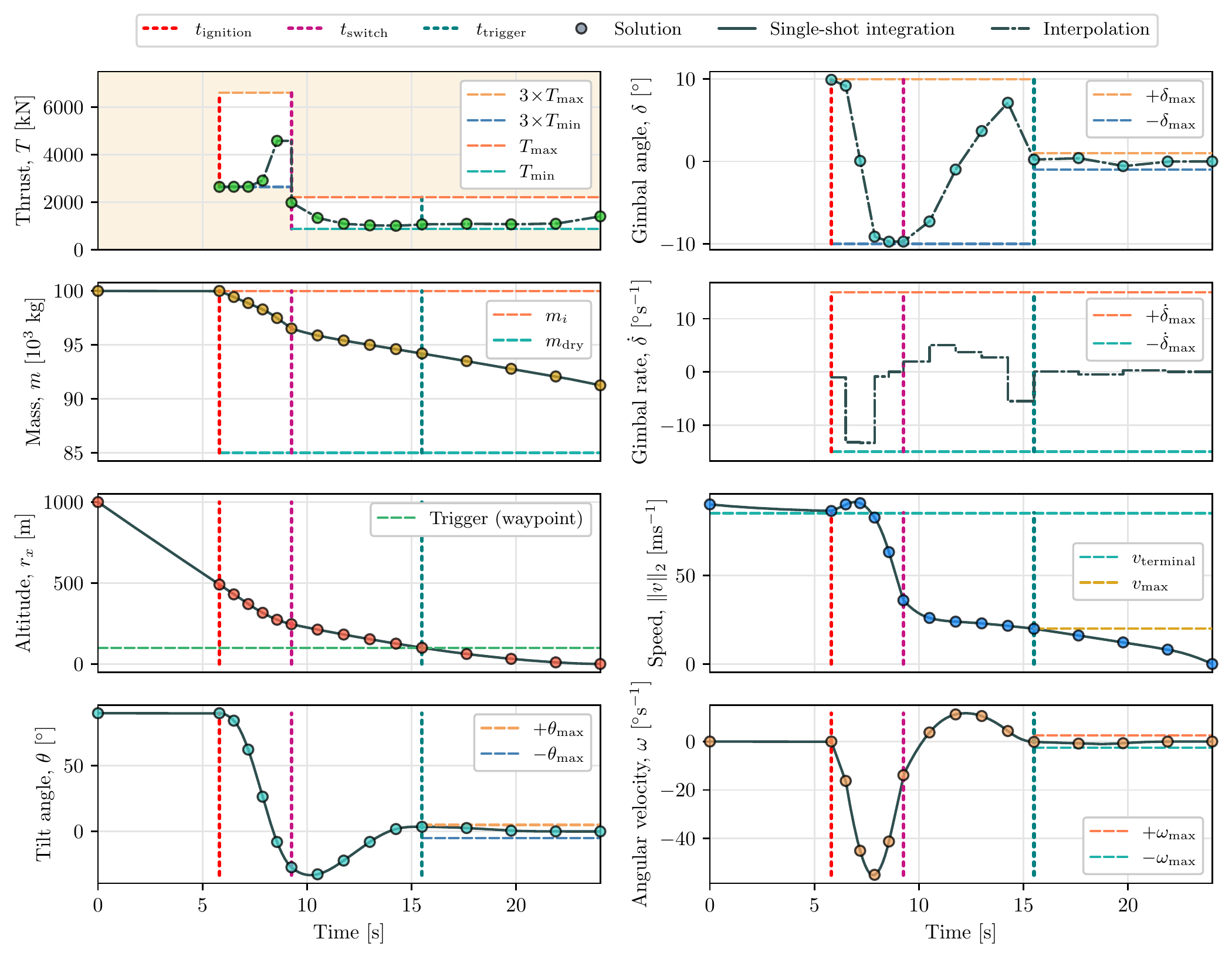}
    \end{minipage}
    \caption{A real-time multi-phase rocket landing guidance solution obtained via {\seco}.}
    \label{fig:starship}
\end{figure*}

\subsection{Common constraints}\label{subsec:common_constraints}
\vspace{-0.375em}
The dynamics and temporal constraints are imposed over the entire horizon. The dynamics constraint is given by Equation \eqref{eq:dyn_constr_disc}. The constraint given by Equation \eqref{eq:grid_constraint} is imposed on the dilation factors to ensure they are bounded.
\begin{align}
    s_{\min} \le s_{l} \le s_{\max},\enskip l \in \{\range{1}{n_{\mathrm{phase}}}\} \label{eq:grid_constraint}
\end{align}
In order to enable the imposition of a single-crossing compound STC in the altitude-triggered terminal descent phase, we impose a minimum altitude constraint in the first three phases, as given by Equation \eqref{eq:min_alt_1}.
\begin{align}
    r^{\xi}_{x_{k}} \ge h_{\mathrm{trigger}},\enskip k \in \{\range{1}{k_\mathrm{trigger}\!-\!1}\} \label{eq:min_alt_1}
\end{align}
\subsection[The unpowered coast phase]{The unpowered coast phase, $k \in \{\range{1}{k_{\mathrm{ignition}}\!-\!1}\}$}\label{subsec:coast}
The center-of-pressure is coincident with the mass-center, i.e., $j=0$ and $l^{j}_{\mathrm{cp}} = l^{0}_{\mathrm{cp}} = 0$ in the dynamics. In implementation, the thrust magnitude is set to zero (the gimbal angle is inconsequential), and a zero-order hold (ZOH) is assumed on the control input signal for this phase, in order to avoid control input constraint violation when the engines are ignited in the next phase. A first-order hold (FOH) can also be assumed here, if the engine startup time is significant and needs to be taken into account (with appropriate constraints on the dilation factor).

The following initial conditions are imposed:
\begin{align}
\begin{split}
    m^{\xi}_{1} = m_{i},\,r^{\xi}_{1} &= r_{i},\,v^{\xi}_{1} = v_{i}, \\
    \theta^{\xi}_{1} = \theta_{i},\,\omega^{\xi}_{1} &= \omega_{i},\,T_{1} = 0 \label{eq:initial_conditions}
\end{split}
\end{align}
where $m_{i}$, $r_{i}$, $v_{i}$, $\theta_{i}$, and $\omega_{i}$ are the initial values for the mass, position, velocity, body tilt angle, and angular velocity, respectively. The thrust magnitude at the first node is constrained to be zero.

\subsection[The high-thrust burn phase]{The high-thrust burn phase, $k \in \{\range{k_{\mathrm{ignition}}}{k_{\mathrm{switch}}\!-\!1}\}$}\label{subsec:high_thrust}
The center-of-pressure shifts towards the nose-cone of the vehicle, i.e., $j=1$ and $l^{j}_{\mathrm{cp}} = l^{1}_{\mathrm{cp}} > 0$ in the dynamics. We model this as a discrete change in its value, which is then held constant for the remainder of the trajectory. The following constraints are imposed on the control variables:
\begin{align}
    &3\,T_{\min} \le T_{k} \le 3\,T_{\max} \label{eq:triple_engine_bounds}\\
\begin{split}
    \max\{-\delta_{\max}, -&\dot{\delta}_{\max} \bar{s}_{k-1} + \bar{\delta}_{k-1}\} \le \delta_{k} \\ &\le \min\{\delta_{\max}, \dot{\delta}_{\max} \bar{s}_{k-1} + \bar{\delta}_{k-1}\} \label{eq:gimbal_combine}
\end{split}
\end{align}
where $T_{\min}$ and $T_{\max}$ are the lower and upper bounds on the thrust magnitude for a single engine, respectively, $-\delta_{\max}$ and $\delta_{\max}$ are the lower and upper bounds on the gimbal deflection angle, respectively, and $-\dot{\delta}_{\max}$ and $\dot{\delta}_{\max}$ are the lower and upper bounds on the gimbal rate, respectively. Equation \eqref{eq:gimbal_combine} combines the gimbal angle and rate constraints, by using reference values in the lower and upper bounds to make them constants and hence, avoid overlapping projections on the same variable. This constraint is, however, exact at convergence, and this formulation has been observed to work well in practice. Although there are only $n_{\mathrm{phase}}$ dilation factor decision variables, we consider the length of $\bar{s}$ to be $N\!-\!1$, such that the dilation factors are repeated to span each phase.

ZOH is assumed on the control input signal between ${k_{\mathrm{switch}}\!-\!1}$ and ${k_{\mathrm{switch}}}$, in preparation for engine downselection, which marks the beginning of the next phase. FOH can also be assumed here, if the engine shutdown time is significant and needs to be accounted for. In order to ensure that the gimbal deflection angle profile does not have any discontinuities, we set the gimbal rate to zero between these nodes, as shown in Equation \eqref{eq:gimbal_hold}.
\begin{align}
\begin{split}
 \delta_{k_{\mathrm{switch}}} = \bar{\delta}_{k_{\mathrm{switch}}-1} \label{eq:gimbal_hold}
\end{split}
\end{align}
\subsection[The low-thrust burn phase]{The low-thrust burn phase, $k \in \{\range{k_{\mathrm{switch}}}{k_{\mathrm{trigger}}\!-\!1}\}$}\label{subsec:low_thrust}
This phase is similar to the high-thrust burn phase, apart from the fact that the switch from 3 to 1 engines occurs at $k_{\mathrm{switch}}$, and the bounds on the thrust magnitude are changed accordingly, as shown in Equation \eqref{eq:single_engine_thrust}. The combined gimbal constraint is left unchanged from the previous phase, and is shown in Equation \eqref{eq:gimbal_combine_unchanged}. FOH is assumed on the control input signal in this phase, and for the remainder of the trajectory.
\begin{align}
    &T_{\min} \le T_{k} \le T_{\max} \label{eq:single_engine_thrust}\\
\begin{split}
    \max\{-\delta_{\max}, -&\dot{\delta}_{\max} \bar{s}_{k-1} + \bar{\delta}_{k-1}\} \le \delta_{k} \\ &\le \min\{\delta_{\max}, \dot{\delta}_{\max} \bar{s}_{k-1} + \bar{\delta}_{k-1}\} \label{eq:gimbal_combine_unchanged}
\end{split}
\end{align}
\subsection[The terminal descent phase]{The terminal descent phase, $k \in \{\range{k_{\mathrm{trigger}}}{N}\}$}\label{subsec:terminal_descent}
The final phase, the terminal descent phase, is the most heavily and tightly constrained phase of flight. This is designed as such in order to enable closed-loop precision landing, i.e, to ensure that the generated guidance trajectories (the feedforward control input signal and the reference state profiles) are amenable to tight tracking via feedback controllers. Such a phase would be especially useful if sub-meter touchdown accuracy is required, for instance, if the vehicle is to be retrieved by the launch tower itself \citep{launchtower}.

The constraint on the thrust magnitude is left unchanged from the previous phase, and is shown in Equation \eqref{eq:single_engine_thrust_td}.
\begin{align}
    &T_{\min} \le T_{k} \le T_{\max} \label{eq:single_engine_thrust_td}
\end{align} 
\begin{subequations}
An altitude-triggered single-crossing compound STC, with five constraint conditions, is imposed, as shown in Equations \eqref{eq:trigger_stc}. The constraint conditions include the following: maximum speed, maximum tilt, maximum angular speed, glideslope, and tighter gimbal deflection bounds. The altitude constraint forms the trigger condition. All of these constraints are imposed in the interval $k \in \{\range{k_{\mathrm{trigger}}}{N\!-1}\}$. 
\begin{align}
    r^{\xi}_{x_{k}} &\left\{
    \begin{aligned} 
        &\!= h_{\mathrm{trigger}},\, \hspace{1.275em}&&\text{if}~k = k_{\mathrm{trigger}} \\
        &\!\le h_{\mathrm{trigger}},\, &&\text{otherwise}
    \end{aligned}\right. \\
    |r^{\xi}_{z_{k}}| \le &\left\{
    \begin{aligned}
        &\tan \gamma_{\mathrm{gs}}\,h_{\mathrm{trigger}}, \hspace{-0.625em}&&\text{if}~k = k_{\mathrm{trigger}} \\
        &\tan \gamma_{\mathrm{gs}}\,\bar{r}^{\xi}_{x_{k}}, &&\text{otherwise}
    \end{aligned} \right.\\
    \|v^{\xi}_{k}\|_{2} &\le v_{\max} \\
    |\theta^{\xi}_{k}| &\le \theta_{\max} \\
    |\omega^{\xi}_{k}| &\le \omega_{\max} \\
    \max\{-\delta_{\max_{\mathrm{TD}}}, -&\dot{\delta}_{\max} \bar{s}_{k-1} + \bar{\delta}_{k-1}\} \le \delta_{k} \\ &\le \min\{\delta_{\max_{\mathrm{TD}}}, \dot{\delta}_{\max} \bar{s}_{k-1} + \bar{\delta}_{k-1}\}\nonumber
\end{align}
\label{eq:trigger_stc}
\end{subequations}
The altitude and glideslope constraints are treated differently at the trigger epoch and after. The altitude constraint is posed as an equality at the trigger—this ensures that the constraint conditions are exactly satisfied at the trigger. Further, the glideslope constraint is cast in the form of box constraints in terms of the reference values, in order to enable closed-form projections (without this measure, there would be two constraints on $r^{\xi}_{x_{k}}$ at every temporal node, thus precluding closed-form projections). Similar to the combined gimbal constraint, this constraint is exact at convergence. The bounds on the gimbal deflection angle are tightened in this phase as well, i.e., $\delta_{\max_{\mathrm{TD}}} < \delta_{\max}$.

The following terminal boundary conditions are imposed:
\begin{align}
\begin{split}
    m^{\xi}_{N} \ge m_{\mathrm{dry}},\,r^{\xi}_{N} &= r_{f},\,v^{\xi}_{N} = v_{f}, \\
    \theta^{\xi}_{N} = \theta_{f},\,\omega^{\xi}_{N} &= \omega_{f},\,\delta_{N} = 0 \label{eq:terminal_conditions}
\end{split}
\end{align}
where $m_{\mathrm{dry}}$ is the dry mass of the vehicle, and $r_{f}$, $v_{f}$, $\theta_{f}$, and $\omega_{f}$ are the terminal values for the position, velocity, body tilt angle, and angular velocity, respectively. The gimbal angle at the final node is constrained to be zero to avoid plume-impingement on the retrieval structure.

\subsection[The discrete SeCO subproblem]{The discrete \begin{normalfont}{\seco}\end{normalfont} subproblem}
We define $J(x_{N})$ in Equation \eqref{eq:objective} to be $-m_{N}$, i.e., the final mass of the vehicle is maximized (thus minimizing propellant consumption). The discrete {\seco} subproblem, which is a second-order cone program (SOCP), can now be given as follows:
\begin{align*}
    \text{min} \quad &\text{\textit{Objective function}: Equation \eqref{eq:objective}}\\
        \text{s.t.} \quad
        & \text{\textit{Common constraints}, Subsection \ref{subsec:common_constraints}:} \nonumber\\
        &\quad \text{Equations \eqref{eq:dyn_constr_disc}, \eqref{eq:grid_constraint}, and \eqref{eq:min_alt_1}}\\
        & \text{\textit{Coast phase}, Subsection \ref{subsec:coast}:} \nonumber\\
        &\quad \text{Equations \eqref{eq:initial_conditions}}\\
        & \text{\textit{High-thrust burn phase}, Subsection \ref{subsec:high_thrust}:} \nonumber\\
        &\quad \text{Equations \eqref{eq:triple_engine_bounds}, \eqref{eq:gimbal_combine}, and \eqref{eq:gimbal_hold}}\\
        & \text{\textit{Low-thrust burn phase}, Subsection \ref{subsec:low_thrust}:} \nonumber\\
        &\quad \text{Equations \eqref{eq:single_engine_thrust} and \eqref{eq:gimbal_combine_unchanged}}\\
        & \text{\textit{Terminal descent phase}, Subsection \ref{subsec:terminal_descent}:} \nonumber\\
        &\quad \text{Equations \eqref{eq:single_engine_thrust_td}, \eqref{eq:trigger_stc}, and \eqref{eq:terminal_conditions}}
\end{align*}
\section{NUMERICAL RESULTS}
\vspace{-0.5em}

For our numerical implementation of the multi-phase rocket landing guidance algorithm, we choose a grid of $N = 16$ discrete temporal nodes, with one node allocated to the unpowered coast phase, and 5 nodes allocated to each of the remaining three phases of flight, i.e, $k_{\mathrm{ignition}} = 2$, $k_{\mathrm{switch}} = 7$, $k_{\mathrm{trigger}} = 12$. The following parameters are used, where the ones pertaining to the vehicle/maneuver were either estimated or obtained from public sources \citep{everydayastronaut, sagliano2021spartan, SCPToolboxCSM2022}.
\begin{gather*}
    g_{0} = 9.81\,\text{m\!\;s$^{-2}$},\,I_{\mathrm{sp}} = 330\,\text{s},\,l_{r} = 4.5\,\text{m},\,l_{h} = 50\,\text{m},\\ l_{\mathrm{cm}} = 0.4\,l_{h},\,l^{0}_{\mathrm{cp}} = 0\,\text{m},\,l^{1}_{\mathrm{cp}} = 0.2\,l_{h},\,\rho_{\mathrm{air}} = 1.225\,\text{kg\!\;m$^{-3}$},\\
    S_{\mathrm{area}} = 545\,\text{m$^2$},\,v_{\mathrm{terminal}} = 85\,\text{m\!\;s$^{-1}$},\,m_{i} = 100000\,\text{kg},\\
    c_{\hat{x}} = 0.0522,\,c_{\hat{z}} = 0.4068,\,T_{\max} = 2200\,\text{kN},\,T_{\min} = 880\,\text{kN},\\
    \delta_{\max} = 10^{\circ},\,\dot{\delta}_{\max} = 15^\circ\text{s$^{-1}$},\,h_{\mathrm{trigger}} = 100\,\text{m},\,\gamma_{\mathrm{gs}} = 5^{\circ},\\
    v_{\max} = 20\,\text{m\!\;s$^{-1}$},\,\theta_{\max} = 5^{\circ},\,\omega_{\max} = 2.5^\circ\text{s$^{-1}$},\,\delta_{\max_{\mathrm{TD}}} = 1^{\circ},\\
    m_{\mathrm{dry}} = 85000\,\text{kg},\,r_{i} = (1000, 100)\,\text{m},\,v_{i} =  
 (-90, 0)\,\text{m\!\;s$^{-1}$},\\
 \theta_{i} = 90^{\circ},\,\omega_{i} = 0^{\circ}\text{s$^{-1}$},\,r_{f} = (0, 0)\,\text{m},\,v_{f} = (0, 0)\,\text{m\!\;s$^{-1}$},\\
 \theta_{f} = 0^{\circ},\,\omega_{f} = 0^{\circ}\text{s$^{-1}$},\,s_{\min} = 0.6\,\text{s},\,s_{\max} = 10\,\text{s}
\end{gather*}
Averaged over 100 full {\seco} solves, we report a mean run-time of the {\pipg} solver (to solve the entire nonconvex problem) of $13.7$ ms. In comparison, {\ecos}\!\footnote{\!\!\!For this problem, {\ecos} is faster than both {\mosek} and {\gurobi}.}\!\! requires $37.1$ ms, on average, to solve the problem. For the comparison, we assess the quality of the converged solutions in terms of the difference in propellant consumption ($0.02\%$) and the number of {\seco} iterations required to converge ($7$). A 
real-time guidance solution obtained via {\pipg} is shown in Figure \ref{fig:starship}. Here, we observe that the optimizer chooses to initiate the powered-descent phase at an altitude of $490.34$ m and a speed of $86.28$ m\!\;s$^{-1}$ (which is very close to the terminal velocity). Further, we note that the numerous prototype flight tests and independent analyses corroborate many of our observations \citep{everydayastronaut}.
\vspace{-0.25em}
\section{CONCLUSIONS}
\vspace{-0.75em}
We introduce {\seco}, a novel matrix-inverse-free paradigm for solving nonconvex optimal control problems in real-time, and solve a multi-phase rocket landing guidance problem with free-transition-time and convex state-triggered constraints. These solutions are computed using {\pipg}, which is almost three times faster than {\ecos}, a state-of-the-art convex optimization solver.
\vspace{-0.125em}


\section*{ACKNOWLEDGEMENTS}
\vspace{-0.75em}
The authors thank the members of the Autonomous Controls Laboratory, especially Taewan Kim, Dayou Luo, and Samet Uzun, and also the members of the Flight Mechanics and Trajectory Design branch (EG5) at the NASA Johnson Space Center, especially Gavin Mendeck, for their valuable insight and the many helpful discussions. The authors also give their special thanks to Miki Szmuk, Taylor Reynolds, and Danylo Malyuta, for laying the foundation for much of this research and also their continued support. This research was supported by NASA grant NNX17AH02A and was partially carried out at the NASA Johnson Space Center; Government sponsorship is acknowledged.
\vspace{-0.25em}


\bibliography{ifacconf}

\end{document}